\documentclass[runningheads,a4paper,11pt]{llncs}
    \usepackage{epsfig}
    \usepackage{amsmath}
    \usepackage{amstext}
    \usepackage{amsfonts}
    \usepackage{amssymb}
\usepackage{algorithmicx}
\usepackage[ruled]{algorithm}
    \usepackage{graphicx}
    \usepackage[margin=1.1in]{geometry}

\newcommand{\ls}[1]
   {\dimen0=\fontdimen6\the\font \lineskip=#1\dimen0
\advance\lineskip.5\fontdimen5\the\font \advance\lineskip-\dimen0
\lineskiplimit=.9\lineskip \baselineskip=\lineskip
\advance\baselineskip\dimen0 \normallineskip\lineskip
\normallineskiplimit\lineskiplimit \normalbaselineskip\baselineskip
\ignorespaces }
\ls{1.0}

\begin{document}
\mainmatter 
\title{Approximately Optimal Scheduling of an M/G/1 Queue with Heavy Tails }
\titlerunning{Truth Serums for Massively Crowdsourced Evaluation Tasks}
%
%
\author{Vijay Kamble%
\and Jean Walrand}
\authorrunning{V. Kamble et al.}
\institute{Department of Electrical Engineering and Computer Sciences\\
 University of California, Berkeley}
%
%
\maketitle
\begin{abstract}
Distributions with a heavy tail are difficult to estimate. If the design of an optimal scheduling policy is sensitive to the details of heavy tail distributions of the service times, an approximately optimal solution is difficult to obtain.  This paper shows that the mean optimal scheduling of an M/G/1 queue with heavy tailed service times does not present this difficulty and that an approximately optimal strategy can be derived by truncating the distributions.

\end{abstract}

\section{Introduction}\label{sec:intro}

Consider the classical stochastic setting of a single server queue with a single service class and the objective of finding a policy that minimizes the mean sojourn time of a job in the system. For this objective, when the service times of the jobs are known at the time of the arrival and when preemption is allowed, the policy of serving the job with the shortest remaining processing time (SRPT) is sample-path optimal irrespective of the distribution of the service times (see \cite{schrage68}). But when the server cannot anticipate the residual service times of the jobs in the system, an age-based optimal policy depends on the service time distribution and is generally difficult to obtain. 

In the case of a G/G/1 queue, if the service time distribution has a decreasing hazard rate (DHR), then the policy of serving the job with the least attained service (LAS) is known to be optimal, while if the distribution satisfies the `new-better-than-used-in-expectation' (NBUE) property\footnote{If $X$ denotes the service time, then the distribution satisfies NBUE if $E(X)>E[X-s\mid X>s]$ for any $s>0$.}, which contains the class of distributions with an increasing hazard rate (IHR), then the optimal policy is the non-preemptive first-come-first-served (FCFS) (see \cite{righter89,righter90,righter92} and also \cite{aar2009}). For classes of distributions between these two extremes, little is known about the optimal policy. In the specific cases of poisson arrivals (i.e. M/G/1) and no arrivals, as a consequence of his seminal result on multi-armed bandit problems, Gittins \cite{gittins} showed that the optimal policy is `index' based. Under this policy, a number that depends on the state of a job, called the Gittins' index, is computed for every job in the system, and the job with the highest index is served at any given time\footnote{An analogous index based optimal policy was derived for the more general case of a multi-class M/G/1 queue, but considering only non-preemptive policies, in the seminal work by Klimov, see \cite{klimov1} and \cite{klimov2}.}. Unfortunately, this index is defined implicitly and it depends on the entire service time distribution of a job. In \cite{aaltoayesta08}, the authors consider the class of distributions that have a NBUE head and a DHR tail (e.g. Pareto distribution), and by analyzing the Gittins' indices, they show that the optimal policy is a combination of the FCFS and LAS service disciplines. Barring these few special cases, there have been no simple characterizations of this optimal policy for the M/G/1 queue\footnote{Recently, \cite{aar2011} gave a remarkable characterization of the Gittins index policy for the queue with no arrivals, where they showed that the optimal policy belongs to the class of multi-level processor sharing disciplines, see \cite{aa2006}.}. 


The sensitivity of the mean optimal policy to the characteristics of the distributions of the service times is undesirable. Even though in principle it is possible to estimate the distribution and calculate the corresponding Gittins' indices, this estimation is difficult for a heavy-tail distribution because of the slow convergence time of the estimates. If the optimal policy is too sensitive to the details of the distribution, inaccuracy in its estimation can potentially have serious consequences. For instance, in the extreme case, although FCFS is optimal for NBUE service distributions, it is known to \emph{maximize} the mean sojourn time for DHR distributions (see again \cite{righter89,righter90,righter92}). Indeed one can show (see \cite{venkat99}) that when the service times have an infinite variance, the mean sojourn time of a job is infinite under any non-preemptive scheduling policy (e.g. FCFS), although simple preemptive policies are able to keep it finite. 

Given this sensitivity, how bad is the consequence of our inability to estimate the tail of the distribution accurately? Can the loss be unbounded? These are the central questions that we try to answer in this paper. Our result is positive: for the specific case of an M/G/1 queue, in which the service times have a finite mean, we show that an approximately optimal policy can be found to the desired magnitude of diminishing error by appropriately truncating the distribution.  Finding the truncation size that corresponds to a desired approximation is simpler than estimating the details of the tail of the distribution of $X$, so this result has practical significance. 

This tendency that policies that perform well for light tailed job service time distributions, perform poorly for heavy tailed distributions and vice versa, has also been observed in the context of a related objective of optimizing the tail of the sojourn time of a job. Here, a scheduling policy needs to asymptotically minimize the probability of facing a long sojourn time and it is assumed that the job sizes are known at arrival. For this objective, on the one hand FCFS is asymptotically optimal for jobs with light tailed service times (see \cite{ramanan01}), where as it performs arbitrarily poorly for heavy tailed jobs (see \cite{venkat99} and \cite{Whitt00}). On the other hand SRPT, LAS, and a recently defined class of policies that appropriately formalize the notion of `prioritizing small jobs' called `SMART' scheduling policies (see \cite{smart1}) perform very well for heavy tailed jobs (i.e. they result in a sojourn time tail that is similar to the service time) but perform poorly if the jobs are light tailed (see \cite{Nuyens06} and \cite{Nuyens08}). A survey of some of these results can be found in \cite{borst03} and  \cite{boxmazwart07}. In fact, it was recently proved in \cite{zwartwierman12} that there can be no scheduling policy in a broadly defined class, which uniformly optimizes the tail of the sojourn time under every distribution of the service time of jobs. Proposing good scheduling policies in the face of this sensitivity is thus important in practice.

\subsection{An overview of the result}
We consider the scheduling of jobs in an M/G/1 queue with a service time distribution $F$ that has a finite mean. At any given time, the server remembers how long he has worked on each job in the queue i.e. the age of the job, and chooses which job to work on. Preemption is allowed. We will focus on work-conserving scheduling policies, which means that as long as there are jobs in the system, the server works at full capacity. The objective of the server is to find a work-conserving, non-anticipative and possibly preemptive scheduling policy that minimizes the mean sojourn time of the jobs in the system. 

Let $X$ be a service time of a typical job and let the rate of arrival be $\lambda$. We show that we can derive an approximately optimal policy from the optimal policy for a related queue with truncated jobs. This truncation is defined as follows.
\begin{definition}\label{typeA}
Let $X$ be the random variable denoting the service time of a typical job with a c.d.f. $F$. Then for any $s\in \mathcal{R}^{+}$, we define a random variable $X^s$ with a c.d.f. $F^s$ defined as
$$
F^s(a)=P(X^s\leq a)=
\left\{
	\begin{array}{ll}
		F(a)  & \mbox{if } a\in [-\infty, s) \\
		1 & \mbox{if } a \geq s
	\end{array}
\right.
$$
This implies that 
$$P(X^s=s)=F^s(s)-\lim_{a\rightarrow s^{-}}F^s(a)=P(X\geq s).$$

\end{definition}
We show that the following policy is approximately optimal. Fix a truncation duration $s$ and compute the optimal policy for scheduling jobs with these truncated service times. Use this policy for scheduling all the jobs that have been served for a duration less than $s$. For any job longer than $s$, after working on it for a duration $s$, place it in a secondary low priority queue. Finally, use the preemptive last-come-first-serve (LCFS) policy to serve this secondary queue as long as there are no jobs left in the primary queue. For this policy, which we call priority optimal-LCFS (PO-LCFS), we show that as the truncation threshold $s$ goes to infinity, the expected contribution to the delay cost from the secondary queue goes to zero. Finding the value of $s$ that corresponds to a desired approximation is easier than estimating the details of the tail of the distribution of $X$. 

Note that the choice of the preemptive LCFS policy to serve the low priority queue was influenced by its ease of analysis. It is possible and quite likely that other preemptive policies (e.g. LAS or processor sharing (PS)) would work equally well\footnote{Using arguments similar to those in \cite{venkat99} it is easy to show that this result would not hold for non-preemptive policies (where preemption is only allowed when a high priority job arrives).}.


\section{Model and Result}
Consider the optimal scheduling problem of an M/G/1 queue in which the service durations of the jobs have a c.d.f. $F$ and the rate of arrival of these jobs is $\lambda$. Since the process describing the total workload in the system is a renewal process that does not depend on the scheduling policy, we consider a single busy period and our aim is to find a work-conserving, non-anticipating, possibly preemptive scheduling strategy $u$, which minimizes the expected value of the sum of sojourn times of all jobs in the busy period. Formally, for a busy period labeled $B$, we want to minimize $E(C(u))$ where
\begin{equation}\label{obj1}
C(u)=\sum_{i=1}^{N_B} S_i(u)
\end{equation}
where $N_B$ is a random variable denoting the number of jobs that arrive in the busy period $B$, and $S_i(u)$ is a random variable denoting the sojourn time of job $i$ when policy $u$ is used to schedule the jobs.
Let $u^*$ be the optimal policy when the service times are distributed according to $F$, and $V^*$ be the corresponding optimal expected cost. We want to define an approximately optimal policy for this scheduling problem, which is derived from the optimal policy for a related problem where the distribution of the job sizes are truncated according to the truncation operation defined in (\ref{typeA}).

Consider those jobs arriving in a busy period $B$ whose service times are larger than $s$. Let $N^s_B$ be the number of these jobs. We can consider each such job to be composed of two smaller jobs, one of length $s$ that arrives when the original job arrives and a second job of the residual length that arrives when the first job of length $s$ is finished. For any fixed scheduling policy, the sojourn time of the original job is exactly the sum of sojourn times of the two smaller jobs. Preserving the label of each original large job for the first small job of service time $s$ that it is composed of and defining new labels for the residual small jobs, we can express the objective in (\ref{obj1}) as minimizing $E(C^s(u))$ where
\begin{equation}\label{obj2}
C^s(u)=\sum_{i=1}^{N_B} S_i(u) +\sum_{j=1}^{N^s_B}S'_j(u).
\end{equation}
The first term of this cost is the contribution by all the  jobs whose service times are smaller than or equal to $s$. The distribution of the service time of each of these jobs is $F^s$. This term can be minimized by using the optimal scheduling policy for the M/G/1 queue in which the service times of the jobs are distributed as $F^s$. In the implementation, the server uses this optimal policy to schedule the jobs and whenever the service time of any job reaches $s$, he stops serving that job and places it in a low priority queue. This event is equivalent to the arrival of the smaller residual job. The low priority queue is thus composed of all the residual jobs and it is served only when there is no job left that has been served for a duration less than $s$. Thus the residual jobs that contribute to the second term in the cost do not interfere with the first term. 
Define the following policy for the original scheduling problem.
\begin{definition} {\bf Priority Optimal + LCFS (PO-LCFS):}
A `priority optimal + LCFS' policy for a truncation parameter $s\in\mathcal{R}^{+}$ 
\begin{enumerate}
\item gives a high priority to jobs that have been served for a duration less than or equal to $s$, while using the optimal scheduling policy for the $M/G/1$ queue in which the service times of the jobs are distributed as $F^s$ for these jobs, and 
\item uses the LCFS policy on jobs that have been served for a duration greater than $s$ only when no high priority jobs remain.
\end{enumerate}
\end{definition}
We will show that under the PO-LCFS policy, denoted by $u^s$, the expected contribution to the total cost in (\ref{obj2}) from the residual jobs diminishes to $0$ as the truncation parameter $s$ grows to infinity. This implies that an approximately optimal policy can the found for the original problem to the desired degree of approximation by choosing an appropriate truncation parameter. Let $V^s$ be the corresponding expected sum of sojourn times of the jobs in a busy period under this policy. 
\begin{theorem}
Let $X$ be the service time of a typical job. Suppose that  $E(X)< \infty$ and $\lambda E(X)<1$. Then under the PO-LCFS policy,
\begin{equation}\label{thm}
\limsup_{s\rightarrow \infty} V^s-V^*\leq \limsup_{s\rightarrow \infty} E(\sum_{j=1}^{N^s_B}S'_j(u^s))=0 
\end{equation}

\end{theorem}
\begin{proof}
First observe that 
\begin{eqnarray}
V^s-V^*&\leq& E\left(\sum_{i=1}^{N_B} S_i(u^s)-\sum_{i=1}^{N_B} S_i(u^*) +\sum_{j=1}^{N^s_B}S'_j(u^s)-\sum_{j=1}^{N^s_B}S'_j(u^*)\right)\nonumber\\
&\leq&E(\sum_{j=1}^{N^s_B}S'_j(u^s)),
\end{eqnarray}
where the inequality follows since the PO-LCFS policy is optimal for jobs with service times smaller than or equal to $s$ and moreover we can ignore the negative term at the end. Now let the sum of sojourn times of all the residual files under the $PO-LCFS$ policy be denoted by $R$ i.e. $R=\sum_{j=1}^{N^s_B}S'_j(u^s)$. We will compute a bound on $R$ using a first step decomposition. To do so, we define a few random variables.

Under the PO-LCFS policy, every time the service duration of a job exceeds $s$, it is placed in a low priority queue and the server continues to work on other smaller jobs (jobs that have been served for a duration less than $s$) in the primary queue. Let $L$ be the first time in the busy period when there are no small jobs left in the primary queue. At $L$, the server starts working on the jobs in the low priority queue for the first time in the busy period. Let $M$ be the number of residual jobs that have accumulated in the low priority queue up until time $L$
and let $Q$ be the total workload, i.e. the sum of the residual service durations, of these $M$ jobs. 

While the server is working on the low priority queue, an arrival of a job in the system causes a digression, since the server switches to work on the new job as dictated by the PO-LCFS policy. This arrival marks the beginning of a small busy period that is stochastically identical to the original busy period. During this busy period, the server does not work on the $M$ jobs that were present in the low priority queue prior to the arrival of this busy period. Once this period ends, the server continues to work on the $M$ jobs until another digression occurs. Let $K(Q)$ be the number of digressions of such complete small busy periods that interrupt the service of the $M$ jobs. Let $W_i$ denote the length of the $i$th small busy period and let $R_i$ be the sum of the sojourn times of all the residual files in the low priority queue under the PO-LCFS policy in the small busy period $i$.
Then $R$ is bounded by the first step decomposition
\begin{equation}
R\leq ML + MQ + M\sum_{i=1}^{K(Q)}W_i +\sum_{i=1}^{K(Q)}R_i.
\end{equation}
The first term in the $R.H.S.$ is the contribution of the delay faced by the $M$ jobs while they wait in the low priority queue before it is served for the first time. This term results from the fact that this delay is no longer than $L$ for each of the jobs. The next term is the contribution by the delay faced by each of the $M$ jobs , while the remaining $M-1$ jobs are being processed. This delay is no larger than the total workload $Q$ of the $M$ jobs, resulting in the second term. Further, each of the $M$ jobs faces a delay due to the digressions, which is no longer than the sum of the durations of all the small busy periods. Finally, the last term is the total delay of all the low priority jobs in the small busy periods. Taking expectations, we get

\begin{equation}\label{mainexp}
E(R)\leq E(ML) + E(MQ) + E(M\sum_{i=1}^{K(Q)}W_i) +E(\sum_{i=1}^{K(Q)}R_i).
\end{equation}
We find upper bounds for each of these terms separately. First, let us consider $E(ML)$. Observe that the random variables $M$ and $L$ do not depend on the particular policy used to schedule the high priority (small) jobs as long as this policy is work conserving. Thus to compute moments of these random variables, we can assume that the preemptive $LCFS$ policy is used for the high priority jobs. Let $X^s$ be the random variable denoting the duration of the first job that arrives in the busy period $B$. Then $M$ can be expressed as
\begin{equation}\label{mexp}
M=\mathbf{1}_{\{X^s=s\}}+\sum_{i=1}^{K(X^s)}M_i,
\end{equation}
where $K(X^s)$ is the number of high priority busy periods that interrupt the service of the first job due to the LCFS policy. These `micro' busy periods are not stochastically identical to the original busy period since during these busy periods, the server only works on the new high priority jobs and returns to serve the original job as soon as there are no high priority jobs left. Thus the duration of this micro busy period $L_i$ has the same distribution as $L$. This also implies that if $M_i$ is the number of residual jobs brought in by the $i$th micro busy period, then $M_i$ has the same distribution as $M$. Thus 
\begin{eqnarray*}
E(M)&=&P(X^s=s)+E(E[\sum_{i=1}^{K(X^s)}M_i\mid X^s])\\
&=& P(X>s)+E(\lambda X^sE(M)).
\end{eqnarray*}
Here the second equality holds because $M_i$ and $X^s$ are mutually independent and conditional on $X^s$, $K(X^s)$ is a poisson random variable with mean $\lambda X^s$. And hence,
\begin{equation}\label{EM}
E(M)=\frac{P(X>s)}{1-\lambda E(X^s)}.
 \end{equation}
Similarly, $L$ can be expressed as 
\begin{equation}
L=X^s+ \sum_{i=1}^{K(X^s)}L_i,
\end{equation}
where $L_i$ is the duration of the $i$th high-priority micro busy period. We can again derive
\begin{equation}\label{EL}
E(L)=\frac{E(X^s)}{1-\lambda E(X^s)},
\end{equation}
which is finite by our assumptions. Further, we have
\begin{eqnarray*}
E(ML)&=&E[(\mathbf{1}_{\{X^s=s\}}+\sum_{i=1}^{K(X^s)}M_i)(X^s+ \sum_{i=1}^{K(X^s)}L_i)]\\
&=&sP(X^s=s)+E[\sum_{i=1}^{K(X^s)}L_i\mid X^s=s]P(X^s=s) \\
&&+~ E(X^s\sum_{i=1}^{K(X^s)}M_i)+ E((\sum_{i=1}^{K(X^s)}M_i)(\sum_{i=1}^{K(X^s)}L_i))\\
&=&sP(X>s)+s\lambda E(L)P(X>s) + E(X^sK(X^s))E(M)\\
&&+ ~E(\sum_{i=1}^{K(X^s)}M_iL_i)+ E(\sum_{i=1}^{K(X^s)}M_i\sum_{j=1;j\neq i}^{K(X^s)}L_j)\\
&=&sP(X>s)+s\lambda E(L)P(X>s) +\lambda E((X^s)^2)E(M)\\
&&+~ \lambda E(X^s)E(ML)+ E(K(X^s)^2-K(X^s))E(M)E(L)\\
&=&sP(X>s)+s\lambda E(L)P(X>s) +\lambda E((X^s)^2)E(M)\\
&&+~ \lambda E(X^s)E(ML)+ \lambda^2E((X^s)^2)E(M)E(L)
\end{eqnarray*}
where the last equality holds because conditional on $X^s$, $K(X^s)$ is a poisson random variable with mean and variance $\lambda X^s$. We thus have 
\begin{equation}\label{EML}
E(ML)= \frac{\left(sP(X>s)+\lambda E((X^s)^2)E(M)\right)(1+\lambda E(L) )}{1-\lambda E(X^s)}
\end{equation}
Let us now consider $E(MQ)$. Notice that we can express $Q$ as,
\begin{equation}\label{qexp}
Q=\sum_{i=1}^{M}\overline{X}_i,
\end{equation}
where $\overline{X}_i$ is the (residual) service time of $i$th job in the low priority queue. Thus we have
\begin{eqnarray}\label{emq}
E(MQ)&=&E(M\sum_{i=1}^{M}\overline{X}_i)\nonumber\\
&=& E(E[M\sum_{i=1}^{M}\overline{X}_i\mid M])\nonumber\\
&=& E(M^2)E(\overline{X}).
\end{eqnarray}
Now from (\ref{mexp}), we can compute $E(M^2)$ as
\begin{eqnarray*}
E(M^2)&=& E[(\mathbf{1}_{\{X^s=s\}}+\sum_{i=1}^{K(X^s)}M_i)^2]\\
&=&P(X>s) +2sP(X>s)\lambda E(M)\\
&&+~ E(\sum_{i=1}^{K(X^s)}M_i^2)+2E(\sum_{1\leq i< j\leq K(X^s)}M_iM_j)\\
&=& P(X>s) +2sP(X>s)\lambda E(M)\\
&&+~ \lambda E(X^s)E(M^2)+E(K(X^s)^2-K(X^s))E(M)^2\\
&=& P(X>s) +2sP(X>s)\lambda E(M)\\
&&+ ~\lambda E(X^s)E(M^2)+\lambda^2E((X^s)^2)E(M)^2.
\end{eqnarray*}
where the last equality again results from the fact that conditional on $X^s$, $K(X^s)$ is a poisson random variable with mean and variance $\lambda X^s$. And thus we have
\begin{eqnarray}\label{EMQ}
E(MQ)&=&E(M^2)E(\overline{X})\nonumber\\
&=& \frac{E(\overline{X})P(X>s)(1+2s\lambda E(M))+E(\overline{X})\lambda^2E((X^s)^2)E(M)^2}{1-\lambda E(X^s)}.
\end{eqnarray}
We will next find a bound for the term $E(M\sum_{i=1}^{K(Q)}W_i)$. We have
\begin{eqnarray}\label{3exp}
E(M\sum_{i=1}^{K(Q)}W_i)&=& E(E[M\sum_{i=1}^{K(Q)}W_i\mid M,Q])\nonumber\\
&=&E(\lambda MQ E(W))\nonumber\\
&=&\lambda E(MQ) \frac{E(X)}{1-\lambda E(X)}
\end{eqnarray}
again because $W_i$ are i.i.d. and conditional on $Q$, $K(Q)$ is a poisson random variable with mean $\lambda Q$. The second equality follows since $E(W)$, which is the expected duration of the busy period $B$, does not depend on the scheduling policy and it is well known, again using first step arguments for the LCFS policy to be $E(W)=\frac{E(X)}{1-\lambda E(X)}$, which is finite by our assumptions. 
Similarly for the last term we have
\begin{eqnarray}\label{4exp}
E(\sum_{i=1}^{K(Q)}R_i)&\leq& E(E[\sum_{i=1}^{K(Q)}R_i)\mid Q])\nonumber\\
&=&E(\lambda Q E(R))\nonumber\\
&=&\lambda E(M) E(\overline{X})E(R).
\end{eqnarray}
Here the second equality follows from (\ref{qexp}).

Thus from (\ref{3exp}) and (\ref{4exp}) and the expression for $E(R)$ in (\ref{mainexp}), we finally have 

\begin{equation}
E(R)\leq \frac{E(ML)+E(MQ)(1+\frac{\lambda E(X)}{1-\lambda E(X)})}{1-\lambda E(M)E(\overline{X})}
\end{equation}
where $E(ML)$ and $E(MQ)$ is bounded by the expressions in (\ref{EML}) and (\ref{EMQ}) respectively. The convergence of this upper bound for increasing threshold $s$ is governed by the following terms.
\begin{enumerate}
\item $E(M)E(\overline{X})=O(P(X>s)E[X-s\mid X>s])=O(P(X>s)E[X\mid X>s]$ from (\ref{EM}).
\item $E(M)E(X^s)=O(P(X>s)(E[X\mid X<s] \,+\, sP(X>s))=O(sP(X>s))$.
\item $E(M)E((X^s)^2)=O(P(X>s)(E(X^2\mathbf{1}_{\{X<s\}}) + s^2P(X>s))$\\ $=O(P(X>s)E(X^2\mathbf{1}_{\{X<s\}})+s^2P(X>s)^2)$.
\end{enumerate}
Now since $E(X)<\infty$ by our assumption, $\lim_{s\rightarrow \infty} P(X>s)E(X\mid X>s)=0$. This also implies that $\lim_{s\rightarrow \infty} sP(X>s)=0$. Further, $E(X^2\mathbf{1}_{\{X<s\}})\leq sE(X \mathbf{1}_{\{X<s\}})\leq E(X)s$ since $E(X)<\infty$ and thus $\lim_{s\rightarrow \infty} P(X>s)E(X^2\mathbf{1}_{\{X<s\}})=0$.
\end{proof}


\section{Conclusion}
We considered the problem of optimal scheduling of an M/G/1 queue with preemption, when the service durations of the jobs are not known at arrival. The goal is to minimize the mean sojourn time of the job in the system. Even though the optimal scheduling policy depends on the entire distribution of the service duration, we proved that an approximately optimal policy can be constructed from an optimal policy for a truncated distribution. Any desired level of approximation error can be achieved by choosing the appropriate truncation size.


\begin{thebibliography}{23}
\footnotesize
\bibitem{aa2006}
Aalto, S., Ayesta, U. (2006). Mean Delay Analysis of Multi Level Processor Sharing Disciplines. In INFOCOM 2006. 25th IEEE International Conference on Computer Communications. Proceedings (pp. 1-11). IEEE.

\bibitem{aaltoayesta08}
Aalto, S., Ayesta, U. (2008). Optimal scheduling of jobs with a DHR tail in the M/G/1 queue. In Proceedings of the 3rd International Conference on Performance Evaluation Methodologies and Tools (p. 50). ICST (Institute for Computer Sciences, Social-Informatics and Telecommunications Engineering).

\bibitem{aar2009}
Aalto, S., Ayesta, U., Righter, R. (2009). On the Gittins index in the M/G/1 queue. Queueing Systems, 63(1-4), 437-458.

\bibitem{aar2011}
Aalto, S., Ayesta, U., Righter, R. (2011). Properties of the Gittins index with application to optimal scheduling. Probability in the Engineering and Informational
 Sciences, 25(03), 269-288.

\bibitem{venkat99}
Anantharam, V. (1999). Scheduling strategies and long-range dependence. Queueing Systems 33, 73-89.
\bibitem{asmussen03}
Asmussen, S. (2003). Applied probability and queues (Vol. 2). New York: Springer.
\bibitem{borst03}
Borst, S. C., Boxma, O. J., Núñez-Queija, R., Zwart, A. P. (2003). The impact of the service discipline on delay asymptotics. Performance Evaluation, 54(2), 175-206.	
\bibitem{boxmazwart07}
Boxma, O., Zwart, B. (2007). Tails in scheduling. ACM SIGMETRICS Performance Evaluation Review, 34(4), 13-20.
\bibitem{gittins}
Gittins, J., Glazebrook, K.,  Weber, R. (1989). Multi-armed Bandit Allocation Indices. John Wiley and Sons, Ltd.

\bibitem{klimov1}
Klimov, G. P. (1975). Time-sharing service systems. I. Theory of Probability and Its Applications, 19(3), 532-551.

\bibitem{klimov2}
Klimov, G. P. (1979). Time-sharing service systems. II. Theory of Probability and Its Applications, 23(2), 314-321.

\bibitem{Nuyens08}
Nuyens, M., Wierman, A., Zwart, B. (2008). Preventing large sojourn times using SMART scheduling. Operations Research, 56(1), 88-101.

\bibitem{Nuyens06}
Nuyens, M., Zwart, B. (2006). A large-deviations analysis of the GI/GI/1 SRPT queue. Queueing Systems, 54(2), 85-97.

\bibitem{righter89}
Righter, R., Shanthikumar, J. G. (1989). Scheduling multiclass single server queueing systems to stochastically maximize the number of successful departures. Probability in the Engineering and Informational Sciences, 3(03), 323-333.

\bibitem{righter90}
Righter, R., Shanthikumar, J. G., Yamazaki, G. (1990). On extremal service disciplines in single-stage queueing systems. Journal of Applied Probability, 409-416.

\bibitem{righter92}
Righter, R., Shanthikumar, J. G. (1992). Extremal properties of the FIFO discipline in queueing networks. Journal of Applied Probability, 967-978.

\bibitem{schrage68}
Schrage, L. (1968). A proof of the optimality of the shortest remaining processing time discipline. Operations Research, 16(3), 687-690.

\bibitem{ramanan01}
Stolyar, A. L., Ramanan, K. (2001). Largest weighted delay first scheduling: Large deviations and optimality. Annals of Applied Probability, 1-48. 
 
\bibitem{Whitt00}
Whitt, W. (2000). The impact of a heavy-tailed service-time distribution upon the M/GI/s waiting-time distribution. Queueing Systems 36, 71-87.
 
\bibitem{smart1}
Wierman, A., Harchol-Balter, M., Osogami, T. (2005). Nearly insensitive bounds on SMART scheduling. In ACM SIGMETRICS Performance Evaluation Review (Vol. 33, No. 1, pp. 205-216). ACM.

\bibitem{zwartwierman12}
Wierman, A., Zwart, B. (2012). Is tail-optimal scheduling possible?. Operations Research, 60(5), 1249-1257.


\end{thebibliography}
\end{document}